\documentclass{article}
\usepackage{latexsym}
\usepackage{amsmath,amsthm,amsfonts,amssymb}

\usepackage{graphicx,latexsym}
\usepackage{lscape}
\usepackage{epstopdf}
\def\keyw{\par\medskip\noindent {\bf Keywords:}\enspace\ignorespaces}
\newtheorem{thm}{Theorem}[section]
 
 \newtheorem{lem}[thm]{Lemma}

\newtheorem{rem}[thm]{Remark}

\numberwithin{equation}{section}

\newtheorem{discu}{Discussion:}

\newtheorem{conje}{Conjecture:}

\begin{document}\date{}
\title{\textbf{Embedding Complete Multipartite Graphs into Certain Trees}}
\author{
A. Arul Shantrinal  $^{a}$
\and
R. Sundara Rajan $^{a,}$\footnote{This work is supported by Project No. 2/48(4)/2016/NBHM-R\&D-II/11580, National Board of Higher Mathematics (NBHM), Department of Atomic Energy(DAE), Government of India.}
\and
 A. Ramesh Babu$^{a}$
\and
S. Anil$^{b}$
\and
Mohammed Ali Ahmed$^{c}$
}
\date{}

\maketitle
\vspace{-0.8 cm}
\begin{center}
$^a$ Department of Mathematics,\\ Hindustan Institute of Technology and Science,\\ Chennai, India, 603 103\\
{\tt shandrinashan@gmail.com }  ~~~~~{\tt vprsundar@gmail.com}  ~~~~~{\tt arbbabu67@gmail.com}\\
\medskip

$^b$ Department of Computer Science and Engineering,\\  Hindustan Institute of Technology and Science,\\ Chennai, India, 603 103\\
{\tt samanthapudianil@gmail.com}\\
\medskip
$^c$ Department of Mathematics, College of Education for Pure Sciences,\\ University of Baghdad, Baghdad, Iraq\\
{\tt mohammedali1975@yahoo.com} \\

\end{center}
\maketitle
\vspace{-0.5 cm}
\normalsize
\begin{abstract}

\noindent One of the important features of an interconnection network is its ability to efficiently simulate programs or parallel algorithms written for other architectures. Such a simulation problem can be mathematically formulated as a graph embedding problem. In this paper, we embed complete multipartite graphs into certain trees, such as $k$-rooted complete binary trees and $k$-rooted sibling trees.

%
\end{abstract}

\vspace{-0.2 cm}
\keyw{Embedding, wirelength, complete multipartite graphs, binary tree, sibling tree}

\vspace{-0.3 cm}
\section{Introduction}
It is well known that a topological structure of an interconnection
network can be modeled by a connected graph $G$. There are a lot
of mutually conflicting requirements in designing the topology of
an interconnection network. It is almost impossible to design a network
which is optimum from all aspects. One of the central issues in designing
and evaluating an interconnection network is to study how well other
existing networks can be embedded into this network and vice-versa. This problem
can be mathematically formulated as a graph embedding problem \cite{XuMa09}.

A graph embedding has two main applications: to transplant parallel
algorithms developed for one network to a different one, and to allocate
concurrent processes to processors in the network. The quality of an embedding can be measured by certain cost criteria. One of these criteria is the \textit{wirelength}. The wirelength of a graph embedding arises from VLSI designs, data structures and data representations, networks for parallel computer systems, biological models that deal with cloning and visual stimuli, parallel architecture, structural engineering and so on {\rm \cite{Xu01,LaWi99}}.

Given two graphs $G$ (guest) and $H$ (host), an embedding from $G$ to $H$ is an injective mapping $f: V(G) \rightarrow V(H)$ and associating a path $P_f(e)$ in $H$ for each edge $e$ of $G$. The wirelength of an embedding $WL(G,H)$ \cite{BeChHaRoSc98} is defined as follows:
$$ WL(G,H)= \underset{f:G\rightarrow H}\min ~~\underset{e=xy\in E(G)}{\sum } \textrm{d}_H(f(x),f(y))=\underset{f:G\rightarrow H}\min ~\underset{e=xy\in E(H)}{\sum } ~EC_{f}(e)$$
where $d_H(f(x),f(y))$ is a distance (need not be a shortest distance) between $f(x)$ and $f(y)$ in $H$ and $EC_{f}(e)$ denote the number of edges $e'$ of $G$
such that $e=xy$ is in the path $P_{f}(e')$ (need not be a shortest path) between $f(x)$ and $f(y)$ in $H$. Further, $EC_{f}(S)=\underset{e\in S}{\sum }EC_{f}(e)$, where $S \subseteq E(H)$.

The \textit{wirelength problem} \cite{BeChHaRoSc98, BeChHaRoSh00,MaRaRaMe09} of a graph $G$ into $H$ is to find an embedding of $G$ into $H$ that induces the minimum wirelength $WL(G,H)$.
The following version of the edge isoperimetric problem of a graph $G(V,E)$ have been considered in the literature \cite{BeDaEl00}, and are $NP$-complete \cite{GaJo79}.

For a subgraph $H$ of $G$ of order $n$,
\vspace{-0.0 cm}
\begin{equation*}
   E_{G}(H) = \{uv \in E(G) ~|~ u, v \in H\}, ~~E_{G}(k)=\underset{H \subseteq  V(G), ~|H|=k} \max ~|E_{G}(H)|
 \end{equation*}

The \textit{maximum subgraph problem} (MSP) for a given $k$, $k\in[n]$ is a problem of computing a
subset $H$ of $V(G)$ such that $\left\vert H\right\vert =k$ and $\left\vert E _{G}(H)\right\vert =E _{G}(k)$. Further, the subsets $H$ are
called the \textit{optimal set} \cite{BeDaEl00, GaJo79, Ha04}.\\

The following results are powerful tools to find wirelength of an embedding using maximum subgraph problem.

\begin{lem}{\rm \cite{MiRaPaRa14}}
\label{modifiedcongestionlemma} Let $f:G\rightarrow H$ be an embedding with $|V(G)|=|V(H)|$. Let $S$ be set of all edges (or edge cut) of $H$ such
that $E(H)\setminus S$ has exactly two components $H_{1}$ and $H_{2}$ and let $G_{i}=[f^{-1}(V(H_{i}))],$ $i=1,2.$ In other words, $G_{i}$ is the induced subgraph on $f^{-1}(V(H_{i}))$ vertices, $i=1,2.$
\textit{Moreover, }$S$ must fulfil the following conditions:
\end{lem}

\vspace{-0.3 cm}
\begin{enumerate}
\item[1.] \textit{For each edge }$ uv\in E(G_{i}), i=1, 2$, $P_{f}(uv)$
\textit{has no edges in the set} $S$.

\item[2.] \textit{For each edge }$ uv\in E(G)$ \textit{with} $u$ \textit{in} $V(G_{1})$ \textit{and} $v$ \textit{in} $V(G_{2})$, $P_{f}(uv)$ \textit{has only one edge in the set} $S$.

\item[3.] $V(G_{1})$ \textit{and} $V(G_2)$ \textit{are optimal sets}.
\end{enumerate}

\noindent \textit{Then }$EC_{f}(S)$\textit{\ is minimum over all embeddings $f:G \rightarrow H$ and }\ $EC_{f}(S)=\underset{v\in V(G_1)}{\sum } deg_G(v)-2|E(G_1)|=\underset{v\in V(G_2)}{\sum } deg_G(v)-2|E(G_2)|$, \textit{where $deg_G(v)$ is the degree of a vertex $v$ in $G$}.

\begin{rem}
\rm{For a regular graph $G$, it is easy to check that, $V(G_2)$ is optimal if $V(G_1)$ is optimal and vice-versa \cite{MaRaRaMe09}.}
\end{rem}

\begin{lem} \label{partitionlemma} {\rm \cite{MiRaPaRa14}}
For an embedding $f$ from $G$ into $H$, let $\{S_{1},S_{2},\ldots,S_{p}\}$ be an edge partition of $[k(H)]$ such
that each $S_{i}$ is an edge cut of $H$ and it satisfies all the conditions of Lemma \ref{modifiedcongestionlemma}. Then%
\begin{equation*}
WL_{f}(G,H)=\frac{1}{k}~\overset{p}{\underset{i=1}{\sum }}EC_{f}(S_{i}).
\end{equation*}
\end{lem}

\section{Main Results}
The multipartite graph is one in all the foremost in style convertible and economical topological structures of interconnection networks. The multipartite has several wonderful options and its one in all the most effective topological structure of parallel processing and computing systems. In parallel computing, a large process is often decomposed into a collection of little sub processes which will execute in parallel with communications among these sub processes. Due to these communication relations among these sub processes the multipartite graph can be applied for avoiding conflicts in the network as well as multipartite networks helps to identify the errors occurring areas in easy way. A complete $p$-partite graph $G = K_{n_1,\ldots, n_p}$ is a graph that contains $p$ independent sets containing $n_i$, $i\in [p]$, vertices, and all possible edges between vertices from different parts.

A tree is a connected graph that contains no cycles. Trees are the most
fundamental graph-theoretic models used in many fields: information theory,
data structure and analysis, artificial intelligence, design of algorithms,
operations research, combinatorial optimization, theory of electrical
networks, and design of networks \cite{RaMaRaAr13}.

The most common type of tree is the binary tree. A binary tree is said to be a
complete binary tree if each internal node has exactly two descendents.
These descendents are described as left and right children of the parent
node. Binary trees are widely used in data structures because they are
easily stored, easily manipulated, and easily retrieved. Also, many
operations such as searching and storing can be easily performed on tree
data structures. Furthermore, binary trees appear in communication pattern
of divide-and-conquer type algorithms, functional and logic programming,
and graph algorithms. A rooted tree represents a data structure with a hierarchical relationship among its various elements \cite{Xu01}.

There are several useful ways in which we can systematically order all nodes
of a tree. Three most important ordering are called
\textit{preorder}, \textit{inorder} and \textit{postorder}. To achieve these
orderings the tree is traversed in a particular fashion. Starting from the
root, the tree is traversed counter clockwise. For preorder, we list a node the first time we pass it.
For inorder, we list a leaf the first time we pass it, but list an interior node the second time we pass it. For postorder, we list a node the last time we pass it  \cite{CoLeRiSt01}.

We now compute the exact wirelength of embedding complete $2^p$-partite graphs $K_{2^{n-p},2^{n-p},\ldots,2^{n-p}}$ into $k$-rooted complete binary trees and $k$-rooted sibling trees for minimizing the wirelength, where $p,n\geq 2$. To prove the main results, we need the following result and the algorithm.

\vspace{0.1 cm}

\begin{lem} \rm{\cite{RaRaLiSe18}}
\label{msplemma}
If $G$ is a complete $p$-partite graph $K_{r,r,\ldots,r}$ of order $pr$, $p,r\geq 2$, then
\setlength{\arraycolsep}{0.001em}
\begin{equation*}
E_G(k) =\left\{
\begin{array}{lcl}
{k(k-1)}/{2} &  \hspace{-3.6cm} ;k\leq p-1 \\\\
{q^2p(p-1)}/{2} &  \hspace{-2.2cm} ;l=qp, ~1\leq q \leq r  \\\\
\frac{(q-1)^2p(p-1)}{2}+j(q-1)(p-1)+\frac{j(j-1)}{2} & \hspace{-2.3cm};l=(q-1)p+j,\\
\hspace{4.7cm}~~~~~~~~~~~~~~ 1\leq j \leq p-1,\\
\hspace{4.57cm}~~~~~~~~~~~~~~ ~2\leq q \leq r.
\end{array}
\setlength{\arraycolsep}{1pt}
\right.
\end{equation*}
\end{lem}

\vspace{-0.3 cm}

\noindent \textbf{Guest Graph Algorithm}\\\\
\textbf{Input}: \hspace{0.1cm} $N = 2^{n}$ (Total number of elements)\\\
\indent \hspace{0.8 cm}  $p \geq 1$, where $2^{n-p}$ represents the number of elements in the each partite\\\\
\textbf{Output}: \hspace{0.1cm} Labeling of complete $2^{p}$-partite graph $K_{2^{n-p},2^{n-p},\ldots,2^{n-p}}$
\small \begin{enumerate}
  \item Begin the algorithm
  \item The guest graph is generated by the complete $2^{p}$-partite graph
  \item The program contains a function $disp_{-}3nr$ which takes \linebreak $2^{n-p}, 2^{n-p},\ldots,2^{n-p}$ as partite elements
  \item Get the values $p$ and $N$, where $p\geq1$
  \item $2^{n-p}$ represents the number of elements in a partite
  \item $2^p$ represents the number of partite generated
  \item \textbf{Assign elements in the partite:}
  \item $m=2^{p}$           \hspace{1.5cm} //Determine number of partite
  \item $y = 0$
  \item for $x \leftarrowР0$ to $n$  do
  \item y ++
  \item $Elem_{-}$val = y
  \item for $i \leftarrowР0$ to $p $ do
  \item for $j \leftarrowР0$ to $p$ do
  \item Array[x]+i+j = $Elem_{-}$val
  \item $Elem_{-}$val = $Elem_{-}$val + $N$
  \item Print the partite:
  \item $n = 0 $  \hspace{1.5cm}                                  //Initiating array number
  \item for $x \leftarrowР0$ to $N$ do
  \item for $i \leftarrow 0$ to $p$ do
  \item	 for $z\leftarrow 0$ to $p$ do
  \item for $j \leftarrow 0$ to $p$ do
 \item Print (array [n] + i + j))
  \item Print a tab space
  \item  n++
  \item  Go to new line
  \item   $z = x \% p$
   \item  if $z = 0$
  \item   Print an empty line
\item 	End the algorithm

\end{enumerate}


\subsection{$k$-rooted Complete Binary Tree}
For any non-negative integer $n$, the complete binary tree of height $n$,
denoted by $T_{n}$, is the binary tree where each internal vertex has
exactly two children and all the leaves are at the same level. Clearly, a
complete binary tree $T_{n}$ has $n$ levels and level $i$, $1\leq i\leq n$,
contains $2^{i-1}$ vertices. Thus $T_{n}$ has exactly $2^{n}-1$ vertices.
The $1$-rooted complete binary tree $T_{n}^{1}$ is obtained from a complete
binary tree $T_{n}$ by attaching to its root a pendant edge. The new vertex
is called the root of $T_{n}^{1}$ and is considered to be at level $0.$ The $%
k$-rooted complete binary tree $T_{n}^{k}$ is obtained by taking $k$ vertex
disjoint $1$-rooted complete binary trees $T_{n}^{1}$ on $2^{n}$ vertices
with roots say $n_{1},n_{2},\ldots,n_{k}$ and adding the edges $(n_{i},n_{i+1})$, $1\leq i\leq k-1$ \cite{RaMaRaAr13}.

%
\paragraph{Wirelength Algorithm A}

\paragraph{Input :}

The complete $2^p$-partite graphs $K_{2^{n-p},2^{n-p},\ldots,2^{n-p}}$, $p,n\geq 2$ and the $1$-rooted complete binary
tree $T_{n}^{1}$ on $2^{n}$ vertices.

\begin{figure}
\centering
\includegraphics[width=7 cm]{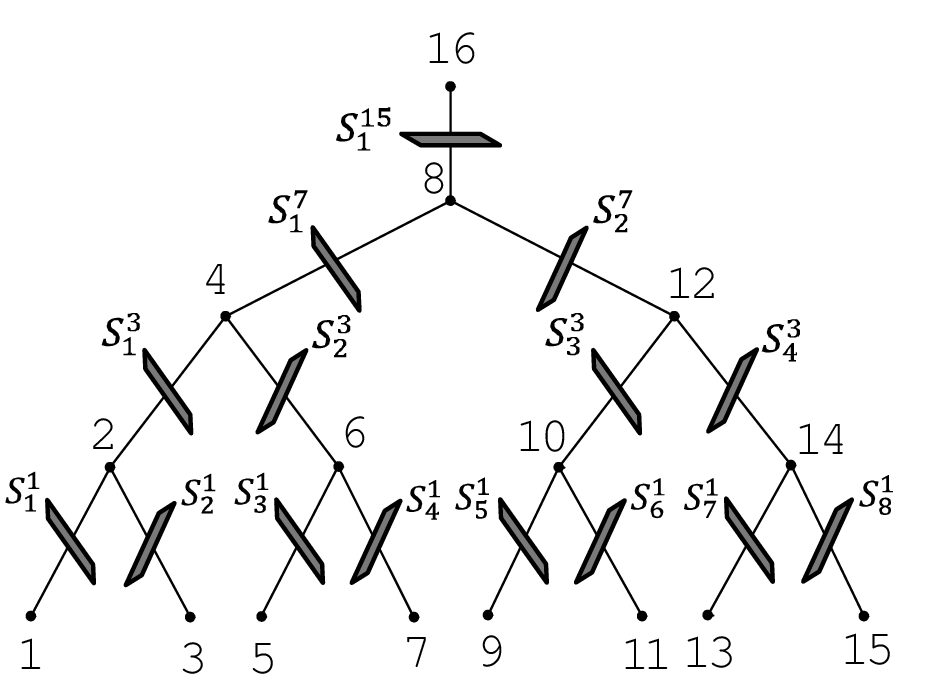}
\caption{Cut edge of 1-rooted complete binary tree $T_4^1$}
\label{fig3}
\end{figure}

\paragraph{Algorithm :}


Label the vertices of complete $2^p$-partite graphs $K_{2^{n-p},2^{n-p},\ldots,2^{n-p}}$ using Guest Graph Algorithm. Label the vertices of $T_{n}^{1}$ using inorder traversal \cite{CoLeRiSt01,RaQuMaWi04} from $1$ to $2^{n}$, see Fig. \ref{fig3}. Let $f(x)=x$ for all $x\in V(K_{2^{n-p},2^{n-p},\ldots,2^{n-p}})$ and for $(a,b)\in E(K_{2^{n-p},2^{n-p},\ldots,2^{n-p}})$, let $P_f(a,b)$ be a shortest path between $f(a)$ and $f(b)$ in $T_{n}^{1}$.

\paragraph{Output :}

An embedding $f$ of $K_{2^{n-p},2^{n-p},\ldots,2^{n-p}}$ into $T_{n}^{1}$ with minimum wirelength and is given by
\begin{eqnarray*}
  WL(G,H) &=&  \left\{
\begin{array}{lcl}
\overset{n}{\underset{j=1}\sum} \overset{2^{n-j}}{\underset{i=1}\sum}(2^j-1)(2^{n-p}(2^p-1)-(2^j-2)) & ; & ~~~j\leq p \\\\
\overset{n}{\underset{j=1}\sum} \overset{2^{n-j}}{\underset{i=1}\sum}(2^p-1)(2^{n-p}(2^j-1)-2^{j-p}(2^j-2)) & ; &~~~j>p
\end{array}
\setlength{\arraycolsep}{1pt}
\right.
\end{eqnarray*}

\paragraph{Proof of correctness :}

For $j=1,2,\ldots,n$ and $i=1,2,\ldots,2^{n-j}$, let $S_{i}^{2^{j}-1}$ be the cut
edge of the $1$-rooted complete binary tree $T_{n}^{1}$, which has one
vertex in level $n-j$ and the other vertex in level $n-j+1$, such that $%
S_{i}^{2^{j}-1}$ disconnects $T_{n}^{1}$ into two components $%
H_{i}^{2^{j}-1} $ and $\overline{H}_{i}^{2^{j}-1}$ where
$V(H_{i}^{2^{j}-1})=\{2^j(i-1)+1,2^j(i-1)+2,\ldots,2^j(i-1)+2^j-1\}$, see Fig. \ref{fig3}. Let $G_{i}^{2^{j}-1} $ and $\overline{G}_{i}^{2^{j}-1}$ be the inverse images of $%
H_{i}^{2^{j}-1}$ and $\overline{H}_{i}^{2^{j}-1}$ under $f$ respectively. By the Guest Graph Algorithm and
Lemma \ref{msplemma}, $G_{i}^{2^{j}-1}$ is an optimal set in $K_{2^{n-p},2^{n-p},\ldots,2^{n-p}}$.
Thus the cut edge $S_{i}^{2^{j}-1}$ satisfies all the conditions of Lemma \ref{modifiedcongestionlemma}. Therefore $EC_{f}(S_{i}^{2^{j}-1})$ is minimum for $j=1,2,\ldots,n$ and $i=1,2,\ldots,2^{n-j}$ and is given by
\begin{equation*}
EC_{f}(S_{i}^{2^{j}-1}) =\left\{
\begin{array}{lcl}
(2^j-1)(2^{n-p}(2^p-1)-(2^j-2)) & ; & ~~~j\leq p \\\\
(2^p-1)(2^{n-p}(2^j-1)-2^{j-p}(2^j-2)) & ; &~~~j>p.
\end{array}
\setlength{\arraycolsep}{1pt}
\right.
\end{equation*}
Then by Lemma \ref{partitionlemma},
\begin{eqnarray*}
  WL(G,H) &=& \overset{n}{\underset{j=1}\sum} \overset{2^{n-j}}{\underset{i=1}\sum}EC_f(S_{i}^{2^{j}-1})\\\\
   &=&  \left\{
\begin{array}{lcl}
\overset{n}{\underset{j=1}\sum} \overset{2^{n-j}}{\underset{i=1}\sum}(2^j-1)(2^{n-p}(2^p-1)-(2^j-2)) & ; & ~~~j\leq p \\\\
\overset{n}{\underset{j=1}\sum} \overset{2^{n-j}}{\underset{i=1}\sum}(2^p-1)(2^{n-p}(2^j-1)-2^{j-p}(2^j-2)) & ; &~~~j>p
\end{array}
\setlength{\arraycolsep}{1pt}
\right.
\end{eqnarray*}

\paragraph{Wirelength Algorithm B}

\begin{figure}
\centering
\includegraphics[width=8 cm]{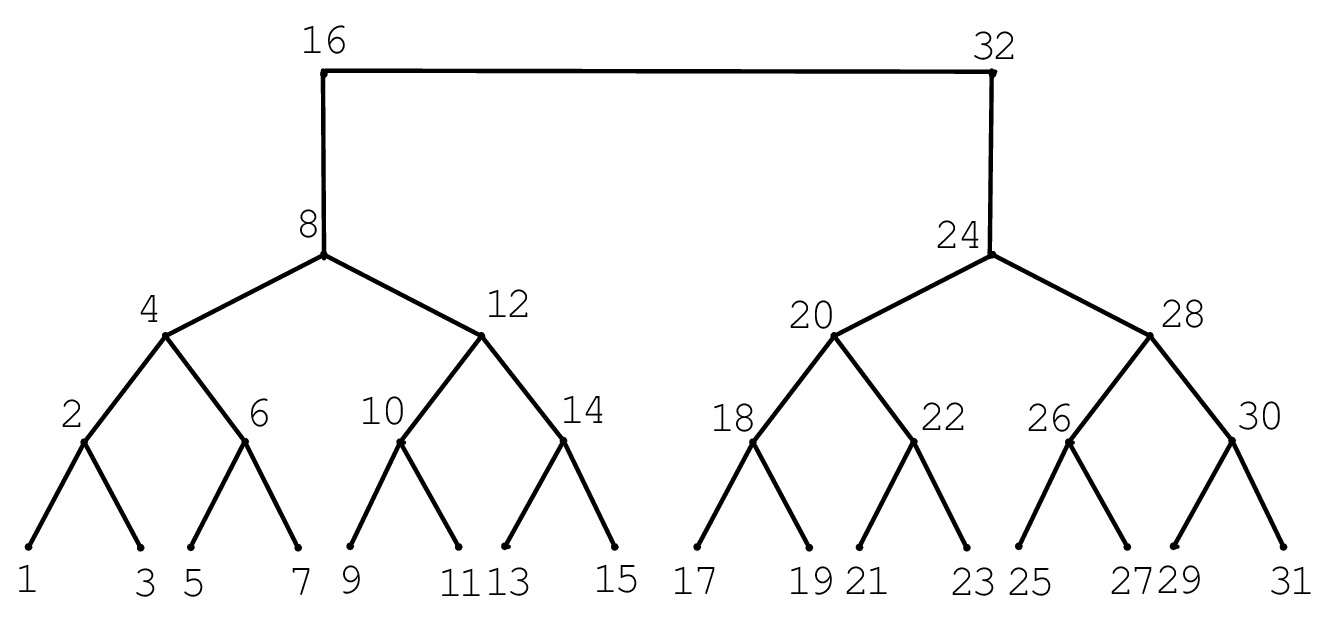}
\caption{Inorder labeling of 2-rooted complete binary tree $T_4^2$}
\label{fig4}
\end{figure}

\paragraph{Input :}

The complete $2^p$-partite graphs $K_{2^{n-p},2^{n-p},\ldots,2^{n-p}}$, $p,n\geq 2$ and the $k$-rooted complete binary
tree $T_{n_1}^{k}$, $k=2^{n-n_1}$.

\paragraph{Algorithm :}

Label the vertices of complete $2^p$-partite graphs $K_{2^{n-p},2^{n-p},\ldots,2^{n-p}}$ using Guest Graph Algorithm. Label the vertices of $T_{n_1}^{k}$, $k=2^{n-n_1}$ as follows:

Let $T_{n_1}^{1,1},T_{n_1}^{1,2},\ldots,T_{n_1}^{1,k}$ be the $k$ vertex disjoint 1-rooted complete binary trees of $T_{n}^{k}$. Label the vertices of $T_{n_1}^{1,i}$, $1\leq i\leq k$ using inorder traversal \cite{CoLeRiSt01,RaQuMaWi04}, see Fig. \ref{fig4}. Let $f(x)=x$ for all $x\in V(K_{2^{n-p},2^{n-p},\ldots,2^{n-p}})$ and for $(a,b)\in E(K_{2^{n-p},2^{n-p},\ldots,2^{n-p}})$, let $P_f(a,b)$ be a shortest path between $f(a)$ and $f(b)$ in $T_{n_1}^{k}$.
\vspace{-0.2cm}
\paragraph{Output :}

An embedding $f$ of $K_{2^{n-p},2^{n-p},\ldots,2^{n-p}}$ into $T_{n_1}^{k}$ with minimum wirelength and is given by
\begin{eqnarray*}
WL(G,H) &=& k\cdot \left\{
\begin{array}{lcl}
\overset{n_1}{\underset{j=1}\sum} \overset{2^{n_1-j}}{\underset{i=1}\sum}(2^j-1)(2^{n-p}(2^p-1)-(2^j-2)) & ; & ~~~j\leq p \\\\
\overset{n_1}{\underset{j=1}\sum} \overset{2^{n_1-j}}{\underset{i=1}\sum}(2^p-1)(2^{n-p}(2^j-1)-2^{j-p}(2^j-2)) & ; &~~~j>p
\end{array}
\setlength{\arraycolsep}{1pt}
\right. \\\\
&&+\left\{
\begin{array}{lcl}
\overset{k-1}{\underset{i=1}\sum}~i2^{n_1}(2^{n-p}(2^p-1)-i2^{n_1}+1) & ; & ~~~i2^{n_1}\leq 2^p \\\\
\overset{k-1}{\underset{i=1}\sum}~i2^{n_1-p}(2^p-1)(2^n-i2^{n_1}) & ; &~~~i2^{n_1}> 2^p
\end{array}
\setlength{\arraycolsep}{1pt}
\right.
\end{eqnarray*}
\paragraph{Proof of correctness :}

By Wirelength Algorithm A, it is enough to prove that the cut edge $(n_{i},n_{i+1})$, $1\leq i\leq k-1 $, where $n_{i}$ is the root of $T_{n_{1}}^{1,i}$, $1\leq i\leq k$, has minimum edge congestion. The cut edge $(n_{i},n_{i+1})$, $1\leq i\leq k-1$
of $T_{n_{1}}^{k}$, disconnects $T_{n_{1}}^{k}$ into two components $H_{i}$
and $\overline{H}_{i}$ where $V(H_{i})=\{1,2,\ldots,i2^{n_{1}}\}$. Let $G_{i}$
and $\overline{G}_{i}$ be the inverse images of $H_{i}$ and $\overline{H}%
_{i} $ under $f$ respectively. By the Guest Graph Algorithm and by Lemma \ref{msplemma}, $G_{i}$ is an
optimal set in $K_{2^{n-p},2^{n-p},\ldots,2^{n-p}}$. Thus the cut edge $(n_{i},n_{i+1})$, $1\leq i\leq
k-1 $ satisfies all the conditions of Lemma \ref{modifiedcongestionlemma}.
Therefore $EC_{f}((n_{i},n_{i+1}))$ is minimum for $i=1,2,\ldots,k-1$ and is given by
\begin{equation*}
EC_{f}((n_{i},n_{i+1})) =\left\{
\begin{array}{lcl}
i2^{n_1}(2^{n-p}(2^p-1)-i2^{n_1}+1) & ; & ~~~i2^{n_1}\leq 2^p \\\\
i2^{n_1-p}(2^p-1)(2^n-i2^{n_1}) & ; &~~~i2^{n_1}> 2^p.
\end{array}
\setlength{\arraycolsep}{1pt}
\right.
\end{equation*}
Then by Lemma \ref{partitionlemma},
\begin{eqnarray*}
  WL(G,H) &=& k\cdot \overset{n_1}{\underset{j=1}\sum} \overset{2^{n_1-j}}{\underset{i=1}\sum}EC_f(S_{i}^{2^{j}-1})+
  \overset{k-1}{\underset{i=1}\sum}EC_f(n_i,n_{i+1})\\\\
    &=& k\cdot \left\{
\begin{array}{lcl}
\overset{n_1}{\underset{j=1}\sum} \overset{2^{n_1-j}}{\underset{i=1}\sum}(2^j-1)(2^{n-p}(2^p-1)-(2^j-2)) & ; & ~~~j\leq p \\\\
\overset{n_1}{\underset{j=1}\sum} \overset{2^{n_1-j}}{\underset{i=1}\sum}(2^p-1)(2^{n-p}(2^j-1)-2^{j-p}(2^j-2)) & ; &~~~j>p
\end{array}
\setlength{\arraycolsep}{1pt}
\right. \\\\
&&+\left\{
\begin{array}{lcl}
\overset{k-1}{\underset{i=1}\sum}~i2^{n_1}(2^{n-p}(2^p-1)-i2^{n_1}+1) & ; & ~~~i2^{n_1}\leq 2^p \\\\
\overset{k-1}{\underset{i=1}\sum}~i2^{n_1-p}(2^p-1)(2^n-i2^{n_1}) & ; &~~~i2^{n_1}> 2^p
\end{array}
\setlength{\arraycolsep}{1pt}
\right.
\end{eqnarray*}
\normalsize

\begin{figure}
\centering
\includegraphics[width=12 cm]{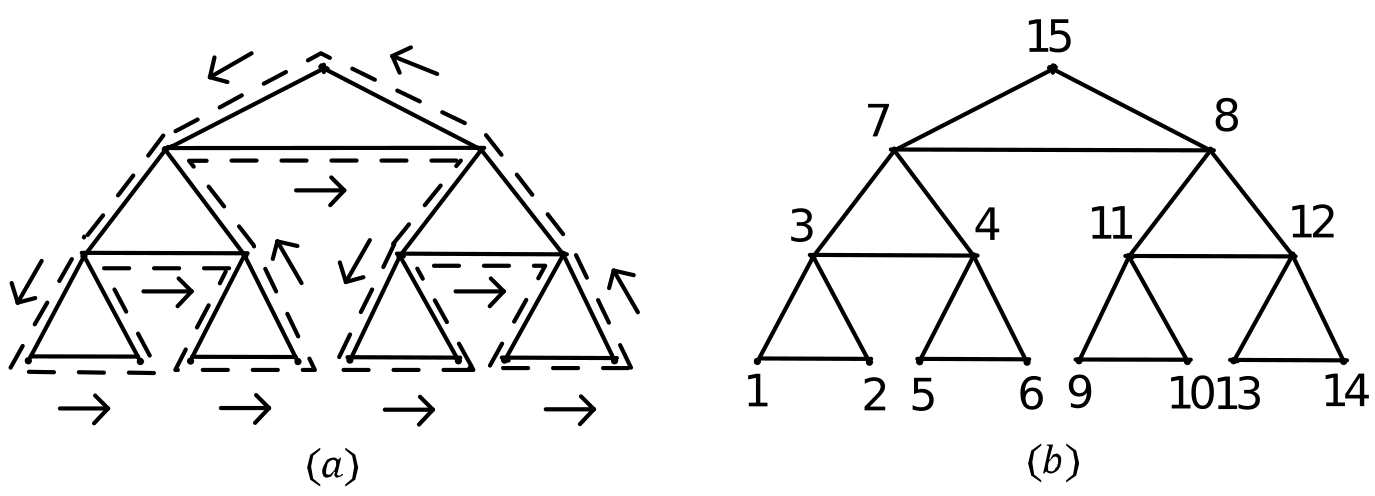}
\caption{(a)~Sibling tree traversal ~~(b)~Labeling of $ST_5$ using sibling tree traversal}
\label{fig5}
\end{figure}

\subsection{$k$-rooted Sibling Tree}

The sibling tree $ST_{n}$ is obtained from the complete binary tree $T_{n}$ by adding edges (sibling edges) between left and right children of the same parent node.
A Sibling tree traversal \cite{RaRaMaMiRa15} follows the usual pattern of binary tree traversal with an additional condition that the traversal does not cut any region, but travels along sibling edges, see Fig. \ref{fig5}(a).


The $1$-rooted sibling tree $ST_{n}^{1}$ is obtained from the $1$-rooted
complete binary tree $T_{n}^{1}$ by adding edges (sibling edges) between
left and right children of the same parent node. The $k$-rooted sibling tree $ST_{n}^{k}$ is obtained by taking $k$ copies of vertex disjoint $1$-rooted sibling tree $ST_{n}^{1}$ on $2^{n}$ vertices
with roots say $n_{1},n_{2},\ldots,n_{k}$ and adding the edges $(n_{i},n_{i+1})$
, $1\leq i\leq k-1$ \cite{RaRaMaMiRa15}.

\paragraph{Wirelength Algorithm C}

\paragraph{Input :}

The complete $2^p$-partite graphs $K_{2^{n-p},2^{n-p},\ldots,2^{n-p}}$, $p,n\geq 2$ and the $1$-rooted sibling tree $ST_{n}^{1}$.

\paragraph{Algorithm :}

Label the vertices of $K_{2^{n-p},2^{n-p},\ldots,2^{n-p}}$ using Guest Graph Algorithm. Label the vertices of $ST_{n}^{1}$ using sibling tree traversal from $1$ to $2^{n}$, see Fig. \ref{fig5}(b). Let $f(x)=x$ for all $x\in V(K_{2^{n-p},2^{n-p},\ldots,2^{n-p}})$ and for $(a,b)\in E(K_{2^{n-p},2^{n-p},\ldots,2^{n-p}})$, let $P_f(a,b)$ be a shortest path between $f(a)$ and $f(b)$ in $ST_{n}^{1}$.

\paragraph{Output :}

An embedding $f$ of $K_{2^{n-p},2^{n-p},\ldots,2^{n-p}}$ into $ST_{n}^{1}$ with minimum wirelength and is given by

\begin{eqnarray*}
 WL(G,H)   &=& \frac{1}{2}\left\{
\begin{array}{lcl}
\overset{n}{\underset{j=1}\sum} ~~\overset{2^{n-j}}{\underset{i=1}\sum}(2^j-1)\\
(2^{n-p}(2^p-1)-(2^j-2)) & ; & ~~~j\leq p \\\\
\overset{n}{\underset{j=1}\sum} ~~\overset{2^{n-j}}{\underset{i=1}\sum}(2^p-1)\\
(2^{n-p}(2^j-1)-2^{j-p}(2^j-2)) & ; &~~~j>p
\end{array}
\setlength{\arraycolsep}{1pt}
\right.
\end{eqnarray*}
\begin{eqnarray*}
\hspace{2.8cm}
 & +&\frac{1}{2}\left\{
\begin{array}{lcl}
\overset{n-1}{\underset{j=1}\sum} ~~\overset{2^{n-j-1}}{\underset{i=1}\sum}(2^{j+1}-2)\\
(2^{n-p}(2^p-1)
-(2^{j+1}-3)) & \hspace{-0.5cm};& j+1\leq p \\\\
\overset{n-1}{\underset{j=1}\sum} ~~\overset{2^{n-j-1}}{\underset{i=1}\sum}(2^p-1)(2^{n-p}(2^{j+1}-2)\\
-2^{j-p+2}(2^j-2))-2 &\hspace{-0.5cm};&  j+1>p.
\end{array}
\setlength{\arraycolsep}{1pt}
\right.\\\\ & +&2^{n-p-1}(2^p-1)
\end{eqnarray*}

\begin{figure}
\centering
\includegraphics[width=6.5 cm]{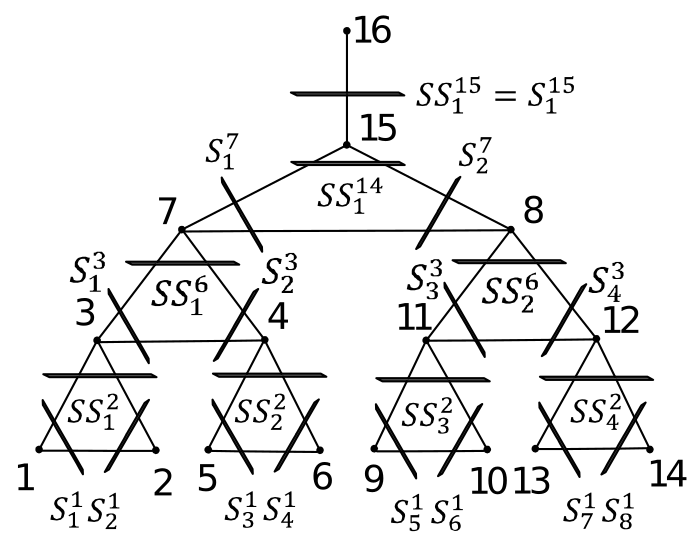}
\caption{Edge cut of 1-rooted sibling tree $ST_4^1$}
\label{fig6}
\end{figure}

\paragraph{Proof of correctness :}
For $j=1,2,\ldots,n$ and $i=1,2,\ldots,2^{n-j}$, let $S_{i}^{2^{j}-1}$ be an edge
cut of the $1$-rooted sibling tree $ST_{n}^{1}$ consisting of edges induced
by the $\left\lceil i/2\right\rceil ^{\text{th}}$ parent vertex from left to
right in level $n-j$ with its left child if $i$ is odd and its right child
if $i$ is even together with the corresponding sibling edge which is the
same edge in either case, such that $S_{i}^{2^{j}-1}$ disconnects $%
ST_{n}^{1} $ into two components $H_{i}^{2^{j}-1}$ and $\overline{H}%
_{i}^{2^{j}-1}$ where $V(H_{i}^{2^{j}-1})$ is consecutively labeled \cite{RaQuMaWi04}, see
Fig. \ref{fig6}. Let $G_{i}^{2^{j}-1}$ and $\overline{G}%
_{i}^{2^{j}-1}$ be the inverse images of $H_{i}^{2^{j}-1}$ and $\overline{H}%
_{i}^{2^{j}-1}$ under $f$ respectively. By Guest Graph Algorithm and Lemma \ref{msplemma},
$G_{i}^{2^{j}-1}$ is an optimal set in $K_{2^{n-p},2^{n-p},\ldots,2^{n-p}}$. Thus the edge cut $%
S_{i}^{2^{j}-1}$ satisfies all the conditions of Lemma \ref{modifiedcongestionlemma}. Therefore $EC_{f}(S_{i}^{2^{j}-1})$ is minimum for $j=1,2,\ldots,n$ and $%
i=1,2,\ldots,2^{n-j}$ and is given by
\begin{equation*}
EC_{f}(S_{i}^{2^{j}-1}) =\left\{
\begin{array}{lcl}
(2^j-1)(2^{n-p}(2^p-1)-(2^j-2)) & ; & ~~~j\leq p \\\\
(2^p-1)(2^{n-p}(2^j-1)-2^{j-p}(2^j-2)) & ; &~~~j>p.
\end{array}
\setlength{\arraycolsep}{1pt}
\right.
\end{equation*}

For $j=1,2,\ldots,n-1$ and $i=1,2,\ldots,2^{n-j-1}$, let $SS_{i}^{2(2^{j}-1)}$ be
an edge cut of the $1$-rooted sibling tree $ST_{n}^{1}$ consisting of the
edges induced by the $i^{\text{th}}$ parent vertex from left to right in
level $n-j$ and its two children, such that $SS_{i}^{2(2^{j}-1)}$
disconnects $ST_{n}^{1}$ into two components $H_{i}^{2(2^{j}-1)}$ and $%
\overline{H}_{i}^{2(2^{j}-1)}$ where $V(H_{i}^{2(2^{j}-1)})$ is consecutively labeled \cite{RaQuMaWi04}, see Fig. \ref{fig6}. Let $%
G_{i}^{2(2^{j}-1)}$ and $\overline{G}_{i}^{2(2^{j}-1)}$ be the inverse
images of $H_{i}^{2(2^{j}-1)}$ and $\overline{H}_{i}^{2(2^{j}-1)}$ under $f$
respectively. By Guest Graph Algorithm and Lemma \ref{msplemma}, $G_{i}^{2(2^{j}-1)}$ is an
optimal set in $K_{2^{n-p},2^{n-p},\ldots,2^{n-p}}$. Thus the edge cut $SS_{i}^{2(2^{j}-1)}$ satisfies
all the conditions of Lemma \ref{modifiedcongestionlemma}. Therefore $EC_{f}(SS_{i}^{2(2^{j}-1)})$ is minimum for $j=1,2,\ldots,n-1$ and $%
i=1,2,\ldots,2^{n-j-1}$. Let $SS_{1}^{2^{n}-1}=S_{1}^{2^{n}-1}$ and it is easy
to see that the conditions of Lemma \ref{modifiedcongestionlemma} are satisfied. We note
that the set $\{S_{i}^{2^{j}-1}:1\leq j\leq n,1\leq i\leq 2^{n-j}\}\cup
\{SS_{i}^{2(2^{j}-1)}:1\leq j\leq n-1,1\leq i\leq 2^{n-j-1}\}\cup
\{SS_{1}^{2^{n}-1}\}$ forms a partition of $[2E(ST_{n}^{1})]$. Further,

\begin{equation*}
EC_{f}(SS_{i}^{2(2^{j}-1)})
 =\left\{
\begin{array}{lcl}
(2^{j+1}-2)\\
(2^{n-p}(2^p-1)-(2^{j+1}-3)) & ; & ~~~j+1\leq p \\\\
(2^p-1)(2^{n-p}(2^{j+1}-2)\\
-2^{j-p+2}(2^j-2))-2 & ; &~~~j+1>p .
\end{array}
\setlength{\arraycolsep}{1pt}
\right.
\end{equation*}

\noindent Then by Lemma \ref{partitionlemma},
\small
\begin{eqnarray*}
  WL(G,H) &=& \frac{1}{2}[\overset{n}{\underset{j=1}\sum} \overset{2^{n-j}}{\underset{i=1}\sum}EC_f(S_{i}^{2^{j}-1})\\
  &+&\overset{n-1}{\underset{j=1}\sum} ~\overset{2^{n-j-1}}{\underset{i=1}\sum}EC_f(SS_{i}^{2(2^{j}-1)})+EC_f(SS_{1}^{2^n-1})]
 \end{eqnarray*}
\begin{eqnarray*}
\hspace{2.3cm}   &=& \frac{1}{2}\left\{
\begin{array}{lcl}
\overset{n}{\underset{j=1}\sum} ~~\overset{2^{n-j}}{\underset{i=1}\sum}(2^j-1)\\
(2^{n-p}(2^p-1)-(2^j-2)) & ; & ~~~j\leq p \\\\
\overset{n}{\underset{j=1}\sum} ~~\overset{2^{n-j}}{\underset{i=1}\sum}(2^p-1)\\
(2^{n-p}(2^j-1)-2^{j-p}(2^j-2)) & ; &~~~j>p
\end{array}
\setlength{\arraycolsep}{1pt}
\right.
\end{eqnarray*}
\begin{eqnarray*}
\hspace{2.8cm}
 & +&\frac{1}{2}\left\{
\begin{array}{lcl}
\overset{n-1}{\underset{j=1}\sum} ~~\overset{2^{n-j-1}}{\underset{i=1}\sum}(2^{j+1}-2)\\
(2^{n-p}(2^p-1)
-(2^{j+1}-3)) & \hspace{-0.5cm};& j+1\leq p \\\\
\overset{n-1}{\underset{j=1}\sum} ~~\overset{2^{n-j-1}}{\underset{i=1}\sum}(2^p-1)(2^{n-p}(2^{j+1}-2)\\
-2^{j-p+2}(2^j-2))-2 &\hspace{-0.5cm};&  j+1>p.
\end{array}
\setlength{\arraycolsep}{1pt}
\right.\\\\ & +&2^{n-p-1}(2^p-1)
\end{eqnarray*}
\normalsize

As $V(ST_{n_1}^k)=V(T_{n_1}^k)$, $k=2^{n-n_1}$, using the proof techniques of Wirelength Algorithm B and Wirelength Algorithm C, we have the following result.

\begin{thm}
The embedding of complete $2^p$-partite graphs \linebreak $K_{2^{n-p},2^{n-p},\ldots,2^{n-p}}$, $p,n\geq 2$ into $k$-rooted sibling tree $ST_{n_1}^k$, $k=2^{n-n_1}$ induces a minimum wirelength $WL(K_{2^{n-p},2^{n-p},\ldots,2^{n-p}},ST_{n_1}^k).$
\end{thm}

\section{Concluding Remarks}
In this paper, we have obtained the wirelength of embedding complete multipartite graphs into certain tree derived architecture, such as $k$-rooted complete binary tree and $k$-rooted sibling trees. Finding the other parameters such as dilation and congestion of embedding complete multipartite graphs into the graphs discussed in this paper are under investigation.


\end{document}